\pdfoutput=1
\documentclass[conference]{IEEEtran}
\IEEEoverridecommandlockouts

\usepackage{amsfonts}
\usepackage{bm}
\usepackage{amsmath,amsfonts,amsthm,amssymb}
\usepackage{setspace}
\usepackage{Tabbing}
\usepackage{fancyhdr}
\usepackage{lastpage}
\usepackage{extramarks}
\usepackage{chngpage}
\usepackage{algorithmic}
\usepackage{soul,color}
\usepackage{epsfig}
\usepackage{graphicx,float,wrapfig,subfigure}
\usepackage{dsfont}
\usepackage{longtable}
\usepackage{wrapfig}
\usepackage{stfloats}
\usepackage{cite}
\usepackage{graphicx}
\usepackage{array}
\usepackage{multirow}{\tiny }
\usepackage{multicol}
\usepackage{colortbl}
\usepackage{tabularx}
\usepackage{mdwmath}
\usepackage{mdwtab}
\usepackage{color}
\usepackage{verbatim}
\usepackage{amsmath}
\usepackage{tikz}
\usetikzlibrary{calc}
\usepackage{flowchart}
\usetikzlibrary{shapes.geometric, arrows}
\usepackage{overpic}
\usepackage{epstopdf}
\usepackage{url}
\usepackage[colorlinks,linkcolor=black,anchorcolor=black,citecolor=black,urlcolor=black]{hyperref} 

\usepackage{makecell}
\usepackage{textcomp}
\usepackage{algorithm}
\usepackage{algorithmic}

\newtheorem{proposition}{Proposition}
\newtheorem{lemma}{Lemma}

\newtheorem{remark}{Remark}

\usepackage{color}
\makeatletter %
\let\myorg@bibitem\bibitem
\def\bibitem#1#2\par{%
	\@ifundefined{bibitem@#1}{%
		\myorg@bibitem{#1}#2\par
	}{%
		\begingroup
		\color{\csname bibitem@#1\endcsname}%
		\myorg@bibitem{#1}#2\par
		\endgroup
	}%
}

\makeatother %

\allowdisplaybreaks

\usepackage{footnote}
\makesavenoteenv{tabular}
\makesavenoteenv{table}

\begin{document}

\title{
	Chance-constrained DC Optimal Power Flow with Non-Gaussian Distributed Uncertainties
\thanks{	
    This paper is funded in part by the Science and Technology Development Fund, Macau SAR (File no. SKL-IOTSC(UM)-2021-2023, File no. 0137/2019/A3 and File no. 0003/2020/AKP). G. Chen's research is also funded in part by the UM Macao PhD Scholarship, University of Macau.
	}
}

\author{
\IEEEauthorblockN{Ge~Chen, Hongcai Zhang, Yonghua Song}
\IEEEauthorblockA{\textit{State Key Laboratory of Internet of Things for Smart City, University of Macau, Macao, China}\\
\textit{hczhang@um.edu.mo}\\
}
}





\maketitle

\begin{abstract}
Chance-constrained programming (CCP) is a promising approach to handle uncertainties in optimal power flow (OPF). However, conventional CCP usually assumes that uncertainties follow Gaussian distributions, which may not match reality. A few papers employed the Gaussian mixture model (GMM) to extend CCP to cases with non-Gaussian uncertainties, but they are only appropriate for cases with uncertainties on the right-hand side but not applicable to DC OPF that containing left-hand side uncertainties. To address this, we develop a tractable GMM-based chance-constrained DC OPF model. In this model, we not only leverage GMM to capture the probability characteristics of non-Gaussian distributed uncertainties, but also develop a linearization technique to reformulate the chance constraints with non-Gaussian distributed uncertainties on the left-hand side into tractable forms. A mathematical proof is further provided to demonstrate that the corresponding reformulation is a safe approximation of the original problem, which guarantees the feasibility of solutions.
\end{abstract}
\begin{IEEEkeywords}
DC optimal power flow, chance-constrained programming, non-Gaussian uncertainties, Gaussian mixture model, linearization
\end{IEEEkeywords}

\section{Introduction}
\IEEEPARstart{I}{n} the past decade, distributed generation (DG) has been widely integrated into power systems to reduce fossil fuel consumption and cut down carbon emission \cite{IEA20202}. However, the growing penetration of DG brings a severe threat to the safety of power systems because its intermittent characteristic may harm the power balance between demand- and generation-side \cite{8550809}. Moreover, since DG is stochastic and hard to predict perfectly, uncertainties will be introduced in optimal power flow (OPF) problem. Therefore, integrating a large amount of DG also exacerbates the difficulty of deciding the power scheduling of generators.

Traditionally, the uncertainties in OPF are handled by robust optimization \cite{8016417}. However, since robust optimization does not allow constraint violation with any realization of uncertainties, the derived solutions are usually overly conservative. Recently, chance-constrained programming (CCP) has become more and more popular \cite{bienstock2014chance,geng2019data}. It only requires constraints that should be satisfied with a pre-determined probability but allows violations in extreme conditions. Thus, the corresponding solution is less conservative but can still guarantee system security in most cases. Reference \cite{8017474} proposed a CCP model to formulate the OPF problem. Reference \cite{pena2020dc} developed a penalty function method to handle the chance constraints in DC OPF. Reference \cite{9535415} built a CCP framework to leverage the flexibility of thermostatically controlled loads to minimize the energy cost of distribution networks. However, most existing works, including \cite{bienstock2014chance,8017474,pena2020dc,9535415}, were based on the Gaussian assumption that uncertainties followed Gaussian distribution. Under this assumption, chance constraints can be reformulated into a second-order cone form. 

In many practical cases, the uncertainties may not be normally distributed, where the previous second-order cone reformulation is no longer valid.
In order to make CCP applicable to non-Gaussian distributed uncertainties, some published works combined a state-of-art statistic technique, Gaussian mixture model (GMM), with CCP in OPF to approximate the original irregular distribution with multiple Gaussian distributions. For instance, reference \cite{7307226} utilized GMM to fit the original irregularly distributed wind generation on the right-hand side and developed a tractable reformulation for the probabilistic OPF. Reference \cite{8936474} also employed GMM to describe the distribution of DG and proposed an analytical reformulation to convert chance constraints into deterministic forms. Reference \cite{8772186} combined GMM with variational Bayesian inference to describe the uncertain wind generation in OPF. In \cite{9376652}, an online-offline double-track distribution fitting approach was proposed to improve the computational efficiency of GMM in probabilistic OPF. In general, these GMM-based methods can better capture the characteristics of uncertainties from DG. They can also reformulate the chance constraints with non-Gaussian distributed uncertainties on the right-hand side into linear forms by finding the corresponding quantile values. However, most of them, including \cite{7307226,8936474,8772186,9376652}, are only suitable for the right-hand side uncertainties. If chance constraints involve non-Gaussian distributed uncertainties on the left-hand side (e.g. chance-constrained DC OPF), this GMM-based method is not applicable (this will be further explained in Section \ref{sec_discuss}). 

The main contribution of this paper is a novel chance-constrained DC OPF model that can handle the non-Gaussian distributed left-hand side uncertainties. This model leverages GMM to capture the probability characteristics of uncertainties. By introducing auxiliary variables, we then develop a bilinear counterpart for the GMM-based chance constraints with left-hand side uncertainties. We further develop a novel linearization technique to reformulate the bilinear terms in the counterpart into linear forms with binary variables that can be directly handled by off-the-shelf solvers. A detailed mathematical proof is also provided to verify that the previous linearized form is a safe approximation of the original model, which guarantees the feasibility of the proposed linearization technique. To the best of our knowledge, this is the first time that the GMM-based chance-constrained method can be extended to handle the DC OPF with left-hand side uncertainties.

The remaining parts are organized as follows. Section \ref{sec_modeling} describes the problem formulation. Section \ref{sec_solution} presents the proposed linearization technique. Section \ref{sec_case} shows simulation results and Section \ref{sec_conclusion} concludes this paper.

\section{Problem formulation} \label{sec_modeling}
\subsection{Chance-constrained DC optimal power flow}
We consider a transmission network with multiple generators and stochastic wind generation.   
Our target is to minimize the total cost of all the generators. By using $i$ to index different buses in the network, the chance-constrained DC optimal power flow problem (CC-OPF) is formulated as follows according to \cite{bienstock2014chance,7973099}:
\begin{align} 
&\min_{\overline{\bm p}, \bm \theta, \bm \alpha} {\mathbb{E}_{\bm \omega} \left(\sum_{i \in \mathcal{G}}c^i(\overline{p}_i - \alpha_{i}(\bm e^\intercal \bm \omega))\right)} \tag{$\textbf{CC-OPF}$},\\
&\begin{array}{r@{\quad}r@{}l@{\quad}l}
\text{s.t.:} &&\sum_{i \in \mathcal{G}} \alpha_{i} = 1, 
\end{array} \label{eqn_alpha_sum}\\
&\sum_{i \in \mathcal{V}} (\overline{p}_i + \mu_i - d_i) = 0, \label{eqn_balance_g_d}\\
& \bm B \bm \theta = \overline{\bm p} + \bm \mu - \bm d, \label{eqn_powerflow}\\
\begin{split}
&\mathbb{P}\left(\beta_{ij}(\overline{\theta}_i-\overline{\theta}_j+[\breve{\bm {B}}(\bm \omega - (\bm e^\intercal \bm \omega)\bm \alpha)]_i \right.\\
&\left. - [\breve{\bm {B}}(\bm \omega - (\bm e^\intercal \bm \omega)\bm \alpha)]_j)\leq f_{ij}^\text{max}\right)\geq 1 - \epsilon,\forall (i,j) \in \mathcal{B},
\end{split} \label{eqn_thermal_limit1}\\
\begin{split}
&\mathbb{P}\left(\beta_{ij}(\overline{\theta}_i-\overline{\theta}_j+[\breve{\bm {B}}(\bm \omega - (\bm e^\intercal \bm \omega)\bm \alpha)]_i) \right.\\
&\left. - [\breve{\bm {B}}(\bm \omega - (\bm e^\intercal \bm \omega)\bm \alpha)]_j)\geq -f_{ij}^\text{max}\right)\geq 1 - \epsilon,\forall (i,j) \in \mathcal{B},
\end{split} \label{eqn_thermal_limit2}\\
&\mathbb{P}\left(\overline{p}_i-(\bm e^\intercal \bm \omega)\alpha_i\geq p_i^\text{min}\right),\forall i \in \mathcal{G}, \label{eqn_power_limit1}\\ &\mathbb{P}\left(\overline{p}_i-(\bm e^\intercal \bm \omega)\alpha_i\leq p_i^\text{max}\right), \forall i \in \mathcal{G}, \label{eqn_power_limit2}
\end{align}
where function $c^i(x)=c_1^i x^2 + c_2^i x$ represents the cost of the generator in bus $i$ with generation $x$. Symbol $\bm \omega \in \mathbb{R}^{|\mathcal{V}|}$ denotes the random fluctuations on nominal wind outputs, where their expectations are zeros (i.e. $\mathbb{E(\bm \omega)}=\bm 0$).
Symbol $\bm B \in \mathbb{R}^{|\mathcal{V}|\times|\mathcal{V}|}$ represents the weighted Laplacian matrix defined as:
\begin{align}
B_{ij} = 
\begin{cases}
-\beta_{ij}, &\text{if } (i,j) \in \mathcal{B},\\
\sum_{k: (k,j) \in \mathcal{B}} \beta_{kj}, &\text{if } i=j, \\
0, &\text{otherwise}. 
\end{cases}
\end{align}
Symbol $\breve{\bm {B}} \in \mathbb{R}^{|\mathcal{V}|\times|\mathcal{V}|}$ is the pseudo-inverse of matrix $\bm B$, which is defined as follows:
\begin{align}
\breve{\bm {B}} = 
\begin{bmatrix} 
\widehat{\bm {B}}^{-1} & 0 \\ 
0 & 0 
\end{bmatrix}
\end{align}
where $\widehat{\bm {B}}$ the submatrix of $\bm B$ by removing the last row and column. Other symbols' definitions can be found in Table \ref{tab_Nomenclature}.
\begin{table}[] 
	\small
	\centering
	\vspace{-4mm}
	\caption{Nomenclature for CC-OPF}
	\begin{tabular}{lll}
		\hline
		\multicolumn{3}{c}{\textbf{Sets}}                                                                    \\ \hline
		$\mathcal{B}$                             & \multicolumn{2}{l}{Set of branches}                                  \\ 
		$\mathcal{G}$                             & \multicolumn{2}{l}{Set of buses with generator}                      \\ 
		$\mathcal{V}$                             & \multicolumn{2}{l}{Set of buses}                                     \\ \hline
		\multicolumn{3}{c}{\textbf{Parameters}}                                                              \\ \hline
		$\beta_{ij}$                      & \multicolumn{2}{l}{susceptance of branch $(i,j)$ ($1/\Omega$)}              \\
		$c_{i}^1,c_{i}^2$                      & \multicolumn{2}{l}{coefficients of generation cost (see Table \ref{tab_gene})}              \\ 
		$d_i$                             & \multicolumn{2}{l}{power demands of bus $i$ (MW)}                               \\ 
		$\bm e$                             & \multicolumn{2}{l}{vector of all ones}                               \\ 
		$f_{ij}^\text{max}$    & \multicolumn{2}{l}{line limit of branch $(i,j)$ (MW)}                                  \\ 
		$\mu_i$                            & \multicolumn{2}{l}{nominal wind power output of bus $i$ (MW)}                  \\ 
		$p_i^\text{max/min}$ & \multicolumn{2}{l}{lower/upper generation bound of bus $i$ (MW)} \\ \hline
		\multicolumn{3}{c}{\textbf{Variables}}                                                               \\ \hline
		$\alpha_i$                         & \multicolumn{2}{l}{proportion of wind power allocated to bus $i$}   \\
		$\overline{p}_i$                             & \multicolumn{2}{l}{nominal output of the generator in bus $i$ (MW)}                     \\ 
		$\theta_i$                         & \multicolumn{2}{l}{nominal voltage angle at bus $i$ (rad)}                     \\ \hline
	\end{tabular} \label{tab_Nomenclature}
	\vspace{-4mm}
\end{table}

Due to the stochastic fluctuations from wind generation, power injections at all buses are uncertain. In order to maintain the equilibrium between generation and demand, {we require that the output of generator $i$ is equal to a nominal value} $\overline{p}_i$ plus the corresponding adjustment $-\alpha_{i}(\bm e^\intercal \bm \omega)$. To ensures that all fluctuations of wind generation are compensated, the summation of $a_i$ should be one, i.e., Eq. (\ref{eqn_alpha_sum}). Constraint (\ref{eqn_balance_g_d}) assures the equilibrium between generation and demand. Eq. (\ref{eqn_powerflow}) represents the DC OPF calculation. Eqs. (\ref{eqn_thermal_limit1})-(\ref{eqn_thermal_limit2}) restrict that the branch power flows are less than the maximum allowable values; Eqs. (\ref{eqn_power_limit1})-(\ref{eqn_power_limit2}) guarantee that the output of each generator keeps in the corresponding feasible region. {The detailed derivation of \textbf{CC-OPF} can be found in  \cite{bienstock2014chance}.}

\subsection{Discussion about CC-OPF} \label{sec_discuss}
All the chance constraints in problem \textbf{CC-OPF} can be written as the following generic form:
\begin{align}
\mathbb{P} \left( \bm a(\bm x)^\intercal \bm \omega \leq b(\bm x)\right) \geq 1-\epsilon. \label{eqn_cc_generic}
\end{align}
Obviously, it contains uncertainties on the left-hand side. If the random vector $\bm \omega$ follows normal distribution, i.e., $\bm \omega\sim \mathcal{N}(\bm \mu, \bm \sigma)$, then  (\ref{eqn_cc_generic}) can be reformulated as:
\begin{align}
\Phi^{-1}(1-\epsilon)\sqrt{\bm a(\bm x)^\intercal \bm \Sigma \bm a(\bm x)} + \bm a(\bm x)^\intercal \bm \mu \leq b(\bm x), \label{eqn_cc_reformulation1}
\end{align}
where $\Phi^{-1}(\cdot)$ is the inverse of the standard normal distribution's cumulative distribution function (CDF).

If  (\ref{eqn_cc_generic}) only contains right-hand side uncertainties, i.e., vector $\bm a$ is not related to $\bm x$, then the original constraint can be directly reformulated as \cite{7307226,8936474,8772186,9376652}:
\begin{align}
\mathcal{Q}_{\bm a^\intercal \bm \omega}^{-1}(1-\epsilon) \leq b(\bm x), \label{eqn_quantile}
\end{align}
where $\mathcal{Q}_{\bm a^\intercal \bm \omega}^{-1}(1-\epsilon)$ represents the quantile of $\bm a^\intercal \bm \omega$ with the violation probability $\epsilon$. Since the quantile has no concern with $\bm x$, it can be calculated based on historical data in advance. 

However, {in practice, chance-constrained DC OPF contains non-Gaussian distributed uncertainties} on the left-hand side, as shown in  (\ref{eqn_cc_generic}). Therefore, either  (\ref{eqn_cc_reformulation1}) or (\ref{eqn_quantile}) can not be directly applied in this task.

\section{Solution Methodology} \label{sec_solution}
\subsection{Gaussian mixture model}
To overcome the aforementioned challenge, we introduce GMM.
The key idea of GMM is to approximate the target distribution with a Gaussian mixture distribution (i.e. a linear combination of multiple Gaussian distributions) based on historical data. Specifically, the PDF of a Gaussian mixture distribution, i.e., $p^\text{GMM}(\bm \xi)$, can be expressed as follows:
\begin{align}
p^\text{GMM}(\bm \omega)=\sum_{j \in \mathcal{J}} \pi_{j}p^\mathcal{N}\left(\bm \zeta_{j}|\bm \mu_{j}, \bm \Sigma_{j}\right), \label{eqn_GMM}
\end{align}
where set $\mathcal{J}$ is the index set of component numbers (i.e., number of Gaussian distributions); $\pi_{j}$ is the weight of component $j$; $p^\mathcal{N}\left(\bm \zeta_{j}|\bm \mu_{j}, \bm \Sigma_{j}\right)$ denotes the PDF of a Gaussian distribution with $\bm \mu_{j}$ as its expectation and $\bm \Sigma_{j}$ as its covariance. Reference \cite{goodfellow2016deep} pointed out that any smooth density can be approximated with any specific, non-zero amount of error by a GMM with enough components.

With historical data, we can use Expectation Maximization method to get the parameters (i.e. weight $\pi_{j}$, expectation $\bm \mu_{j}$, and covariance $\bm \Sigma_{j}$) in  (\ref{eqn_GMM}). To obtain the best fitting results, a cross-validation technique is also adopted to tune the hyper-parameter $J$. Then, the target distribution, i.e. actual distribution of $\bm \omega$, can be well approximated with  (\ref{eqn_GMM}).

\subsection{Counterpart of GMM-based chance constraints}
We decompose the original chance constraint (\ref{eqn_cc_generic}) into a linear combination of multiple chance constraints with normally distributed uncertainties. This decomposition is based on the following Lemma \cite{hu2018chance}. 

\begin{lemma} \label{lemma_1}
If the distribution of the random variable $\bm \omega$ can be expressed as a Gaussian mixture distribution, i.e., Eq. (\ref{eqn_GMM}), then we have 
\begin{align}
\mathbb{P} \left( \bm a(\bm x)^\intercal \bm \omega \leq b(\bm x)\right)=\sum_{j \in \mathcal{J}}\pi_j\mathbb{P} \left( \bm a(\bm x)^\intercal \bm \zeta_{j} \leq b(\bm x)\right).
\end{align}
\end{lemma}

Based on \textbf{Lemma} \ref{lemma_1}, Eq. (\ref{eqn_cc_generic}) can be reformulated as:
\begin{align}
\sum_{j \in \mathcal{J}}\pi_{j} \mathbb{P}\left( \bm a(\bm x)^\intercal \bm \zeta_{j} \leq b(\bm x)\right) \geq 1 - \epsilon. \label{eqn_cc_reformulation2}
\end{align}
By introducing auxiliary variables $y_{j}$ for each component of the GMM, we can get the counterpart of  (\ref{eqn_cc_reformulation2}), as follows:
\begin{align}
&\sum_{j \in \mathcal{J}}\pi_{j}y_{j} \geq 1-\epsilon, \quad y_{j}\geq 0.5,\quad \forall j \in \mathcal{J}, \label{eqn_reformulation_y} \\
&\Phi^{-1}(y_{j})\sqrt{\bm a(\bm x)^\intercal \bm \Sigma_j \bm a(\bm x)} + \bm \mu_j^\intercal a(\bm x) \leq b(\bm x), \forall j \in \mathcal{J}. \label{eqn_reformulation_3}
\end{align}
where $\Phi^{-1}(\cdot)$ is the inverse of the standard normal distribution's CDF. Eq. (\ref{eqn_reformulation_3}) eliminates the intractable probability. However, it is still nonconvex for off-the-shelf solvers due to term $\Phi^{-1}(y_{j})\sqrt{\bm a(\bm x)^\intercal \bm \Sigma_j \bm a(\bm x)}$ (note $y_{j}$ is also a variable). Moreover, the explicit mathematical expression of $\Phi^{-1}(\cdot)$ is also hard to obtain, which further increases the difficulty of the power dispatch. 

\subsection{Linearization technique for intractable terms} \label{sec_linearization}
To make the counterpart (\ref{eqn_reformulation_3}) tractable, we propose a novel linearization technique to convert the intractable terms into linear and second-order cone constraints with binary variables. Then, the OPF problem becomes a mixed-integer program that can be efficiently solved by off-the-shelf solvers.

Firstly, we introduce two auxiliary variable $s_{j}$ and $t_{j}$ to re-express  (\ref{eqn_reformulation_3}) into the following form:
\begin{align}
&\Phi(s_{j}) \geq y_j, \quad \forall j \in \mathcal{J}, \label{eqn_inverse}\\
&\sqrt{\bm a(\bm x)^\intercal \bm \Sigma_j \bm a(\bm x)} \leq t_j, \quad \forall j \in \mathcal{J},\label{eqn_reformulation_4_socp}\\
&s_j \cdot t_j + \bm \mu_j^\intercal \bm a(\bm x) \leq b(\bm x), \quad \forall j \in \mathcal{J}. \label{eqn_reformulation_4_bilinear}
\end{align}
Obviously, the above inequalities are an inner approximation of  (\ref{eqn_reformulation_3}). Among the above three constraints, Eq. (\ref{eqn_reformulation_4_socp}) is convex (second-order cone constraint). However, Eq. (\ref{eqn_inverse}) is intractable since the CDF $\Phi(s_{j})$ is not an elementary function, i.e., it contains integration operator, which can not be handled by off-the-shelf solvers. Moreover, Eq. (\ref{eqn_reformulation_4_bilinear}) is non-convex due to the bilinear term $s_j \cdot t_j$. 

\begin{figure}
	\centering
	\subfigbottomskip=-6pt
	\subfigcapskip=-4pt
	\subfigure[]{\includegraphics[width=0.49\columnwidth]{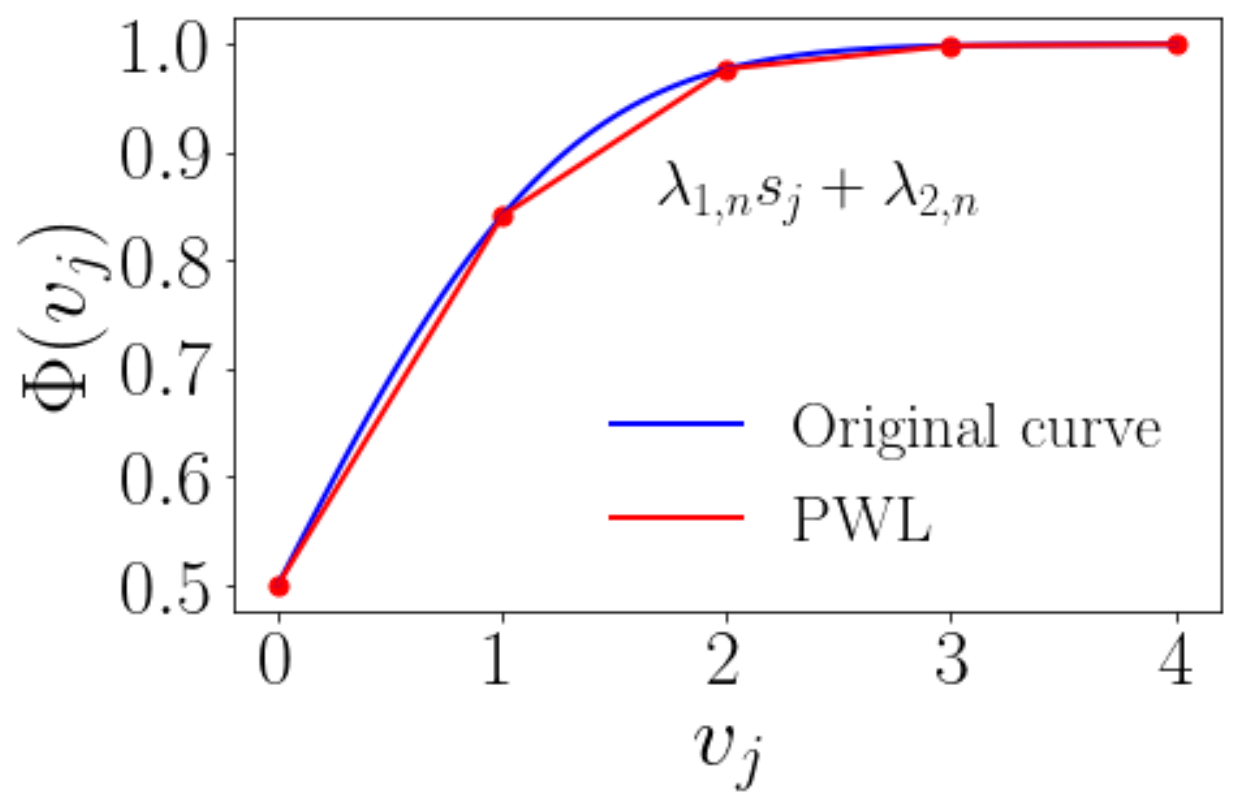}}
	\subfigure[]{\includegraphics[width=0.47\columnwidth]{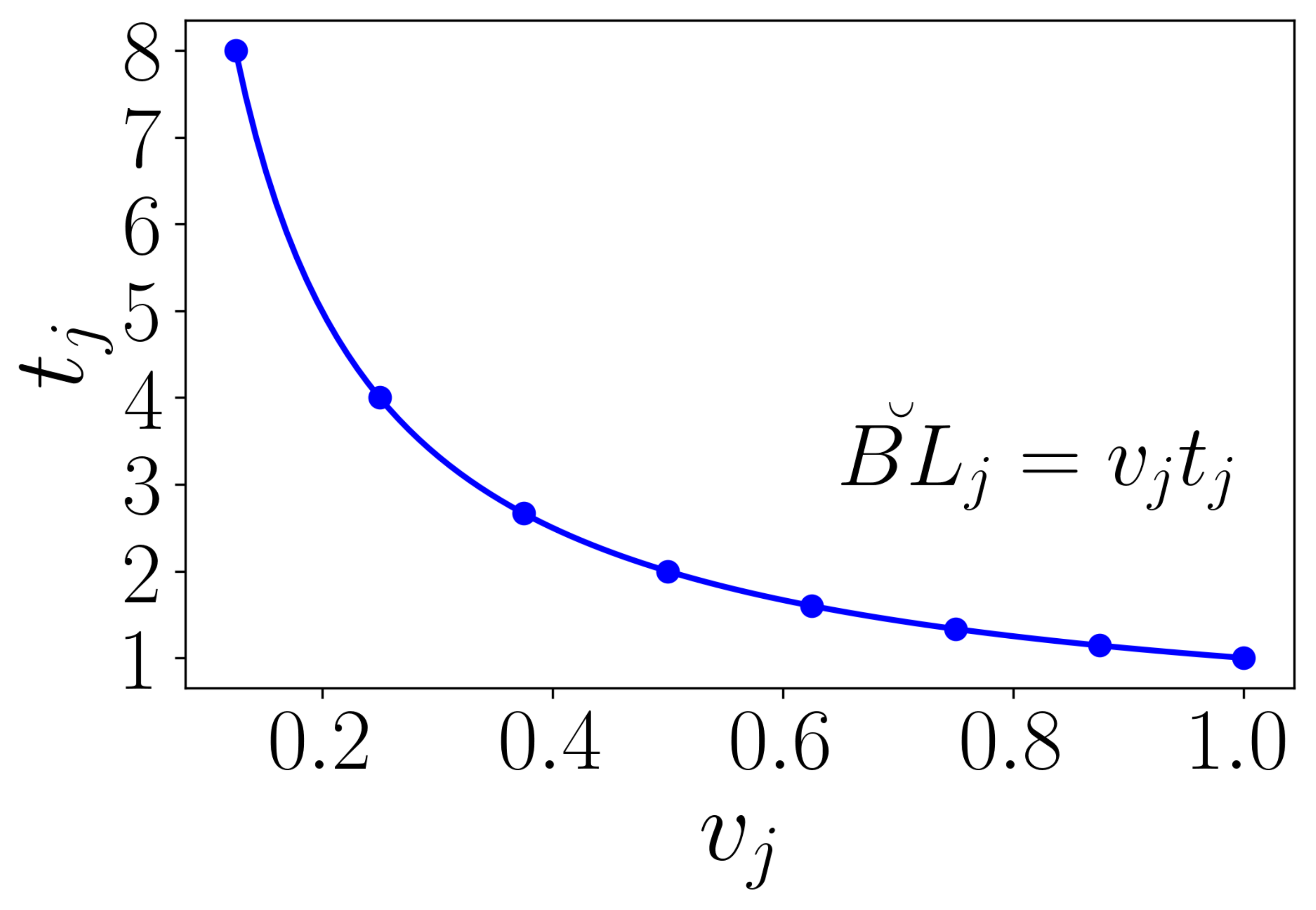}}
	\caption{Schematic diagram of (a) the piece-wise linearization on $\Phi(s_j)$ and (b) discretization on $v_j$. In (a), function $\lambda_{1,n}s_j+\lambda_{2,n}$ represents the $n$-th line segment after piece-wise linearization. Obviously, function $\Phi(s_{j})$ is concave. PWL means piece-wise linearization. In (b), we use the discrete points to replace the original feasible set (blue line).
	}
	\label{fig_ppf}
\end{figure}

To make  (\ref{eqn_inverse}) tractable, we adopt piece-wise linearization to approximately represent $\Phi(s_{j})$. Then, Eq. (\ref{eqn_inverse}) can be approximated by:
\begin{align}
\min\{\lambda_{1,n}s_j + \lambda_{2,n}, \forall n \in \mathcal{N}\} \geq y_j, \quad \forall j \in \mathcal{J}, \label{eqn_inverse2}
\end{align}
where $\lambda_{1,n}$ and $\lambda_{2,n}$ define the $n$-th line segment, which is constructed by connecting two different points on the function $\Phi(s_{j})$, as shown in Fig. \ref{fig_ppf}(a). Set $\mathcal{N}=\{1,2,\cdots,N\}$ is the line segment number. Note  (\ref{eqn_inverse2}) is a linear constraint and tractable for off-the-shelf solvers.

\begin{proposition}
	Eq. (\ref{eqn_inverse2}) is an inner approximation of  (\ref{eqn_inverse}).
\end{proposition} 
\emph{Proof}: The Hessian matrix of function $\Phi(s_{j})$ is
\begin{align}
\Phi^{''}(s_{j}) = \frac{1}{\sqrt{2\pi}}\cdot(-s_j)\cdot e^{-s_j^2/2} \leq 0, \forall s_j \in \mathbb{R}_+.
\end{align}
Thus, when $s_j\geq 0$ ($y_j \geq 0.5$), function $\Phi(s_{j})$ is concave. According to the definition of the concave function, we have
\begin{align}
\psi \Phi(x_1) + (1-\psi)\Phi(x_2) \leq \Phi(\psi x_1 + (1-\psi) x_2).
\end{align}
That is to say, the graph of $\Phi(s_{j})$ is always above the previous line segments, so 
\begin{align}
\Phi(s_{j}) \geq \min\{\lambda_{1,n}s_j + \lambda_{2,n}, \forall n \in \mathcal{N}\}, \forall s_j \in \mathbb{R}_+.
\end{align}

We further develop a safe approximation for  (\ref{eqn_reformulation_4_bilinear}). The key idea is to discretize one variable $s_j$ and then employ the big-M method to restrict the value of the bilinear term. Supposing that the lower and upper bounds of $s_j$ are $s^\text{max}$ and $s^\text{min}$, respectively, variable $s_j$ can be standardized with a new variable $v_j \in [0,1]$:
\begin{align}
s_j = (s^\text{max}-s^\text{min})v_j + s^\text{min}. \label{eqn_stand}
\end{align}
The value of $v_j$ can be expressed as a discrete form:
\begin{align}
v_j = \sum_{l \in \mathcal{L}} 2^l \cdot z_{j,l}, \quad  z_{j,l} \in \{0,1\}, \quad \forall l \in \mathcal{L}, \label{eqn_discrete}
\end{align}
where $\mathcal{L}=\{LB,LB+1,\cdots,UB\}$ is the index set of powers; parameters $LB=-\infty$ and $UB=-1$. By using $\breve{BL}_j$ to denote the bilinear term $v_jt_j$, we have
\begin{align}
&\breve{BL}_j = \sum_{l \in \mathcal{L}} 2^l \cdot \breve{bl}_{j,l}, \label{eqn_bilinear} \\
&t_\text{min} z_{j,l} \leq \breve{bl}_{j,l} \leq t_\text{max} z_{j,l}, \forall l \in \mathcal{L},\\
& - t_\text{max} (1-z_{j,l}) \leq \breve{bl}_{j,l}-t_j \leq t_\text{max} (1-z_{j,l}), \forall l \in \mathcal{L}.
\end{align}
Then, the bilinear constraint (\ref{eqn_reformulation_4_bilinear}) can be reformulated as
\begin{align}
(s^\text{max}-s^\text{min})\breve{BL}_j+s^\text{min}t_j + \bm \mu_j^\intercal \bm a(\bm x) \leq b(\bm x),\forall j \in \mathcal{J}. \label{eqn_bilinear_reformulation}
\end{align}
Essentially, the proposed discretization method chooses multiple discrete points to replace the original continuous feasible region of $s_j$, as shown in Fig. \ref{fig_ppf}(b). Note since we shrink the feasible region of $s_j$, this discretization is also an inner approximation method.
The stepsize for choosing points is $2^{LB}$. If $LB=-\infty$, the set formed by the chosen discrete points is the same with the original feasible region. Nevertheless, with a relatively large stepsize (e.g. $LB=-4$), guaranteed optimality can be also achieved.


\subsection{Reformulation of CC-OPF}
Observing that the expectation of uncertain parameter $\bm \omega$ is zero, the objective is equivalent to
\begin{align}
&\mathbb{E}_{\bm \omega} \left(\sum_{i \in \mathcal{G}}c^i(\overline{p}_i - \alpha_{i}(\bm e^\intercal \bm \omega))\right) = \sum_{i \in \mathcal{G}} \mathbb{E}_{\bm \omega} \left(c^i(\overline{p}_i - \alpha_{i}(\bm e^\intercal \bm \omega))\right) \notag\\
& = \sum_{i \in \mathcal{G}} \left(c^i_1 \overline{p}_i^2 + c^i_2 \overline{p}_i +c^i_1 \alpha_{i}^2 \bm e^\intercal \bm \Sigma \bm e \right),
\end{align}
where $\bm \Sigma$ is the covariance matrix of $\bm \omega$. Based on the linearization technique for  (\ref{eqn_inverse}) and discretization method for the counterpart (\ref{eqn_reformulation_4_bilinear}), all the intractable terms can be removed. Then, the \textbf{CC-OPF} can be reformulated as:
\begin{align} 
&\min_{\overline{\bm p}, \bm \theta, \bm \alpha} {\sum_{i \in \mathcal{G}} \left(c^i_1 \overline{p}_i^2 + c^i_2 \overline{p}_i +c^i_1 \alpha_{i}^2 \bm e^\intercal \bm \Sigma \bm e \right)} \tag{$\textbf{CC-OPF-2}$},\\
&\begin{array}{r@{\quad}r@{}l@{\quad}l}
\text{s.t.:} &&\text{Eqs. (\ref{eqn_alpha_sum}), (\ref{eqn_powerflow}), \{(\ref{eqn_reformulation_4_socp}), (\ref{eqn_inverse2}), (\ref{eqn_discrete})-(\ref{eqn_bilinear_reformulation})\}$_{\forall m \in \mathcal{M}}$.} 
\end{array} \notag
\end{align}
Note in  \textbf{CC-OPF} we have multiple chance constraints, i.e., Eqs. (\ref{eqn_thermal_limit1})-(\ref{eqn_power_limit2}). Each chance constraint should be linearized into  (\ref{eqn_reformulation_4_socp}), (\ref{eqn_inverse2}), (\ref{eqn_discrete})-(\ref{eqn_bilinear_reformulation}) based on the technique proposed in Section \ref{sec_linearization}. 
Thus, we employ $m$ to index the original chance constraint and use   \{(\ref{eqn_reformulation_4_socp}), (\ref{eqn_inverse2}), (\ref{eqn_discrete})-(\ref{eqn_bilinear_reformulation})\}$_{\forall m \in \mathcal{M}}$ to cover reformulations for all chance constraints. 

\begin{remark}
	\textbf{CC-OPF-2} is a safe approximation of \textbf{CC-OPF} because both the linearization technique for the CDF $\Phi(\cdot)$ and discretization method for bilinear term $s_j t_j$ are inner approximations. Thus, a feasible solution of \textbf{CC-OPF-2} must be also feasible for \textbf{CC-OPF}.
\end{remark}

\begin{remark}
	\textbf{CC-OPF-2} is a mixed-integer second-order cone program. It can be efficiently solved by some mature algorithms (e.g. branch-and-bound) in off-the-shelf solvers.
\end{remark}

\section{Case study} \label{sec_case}
\subsection{Simulation setting up}
\begin{figure}
	\centering
	\vspace{-4mm}
	\includegraphics[width=0.7\columnwidth]{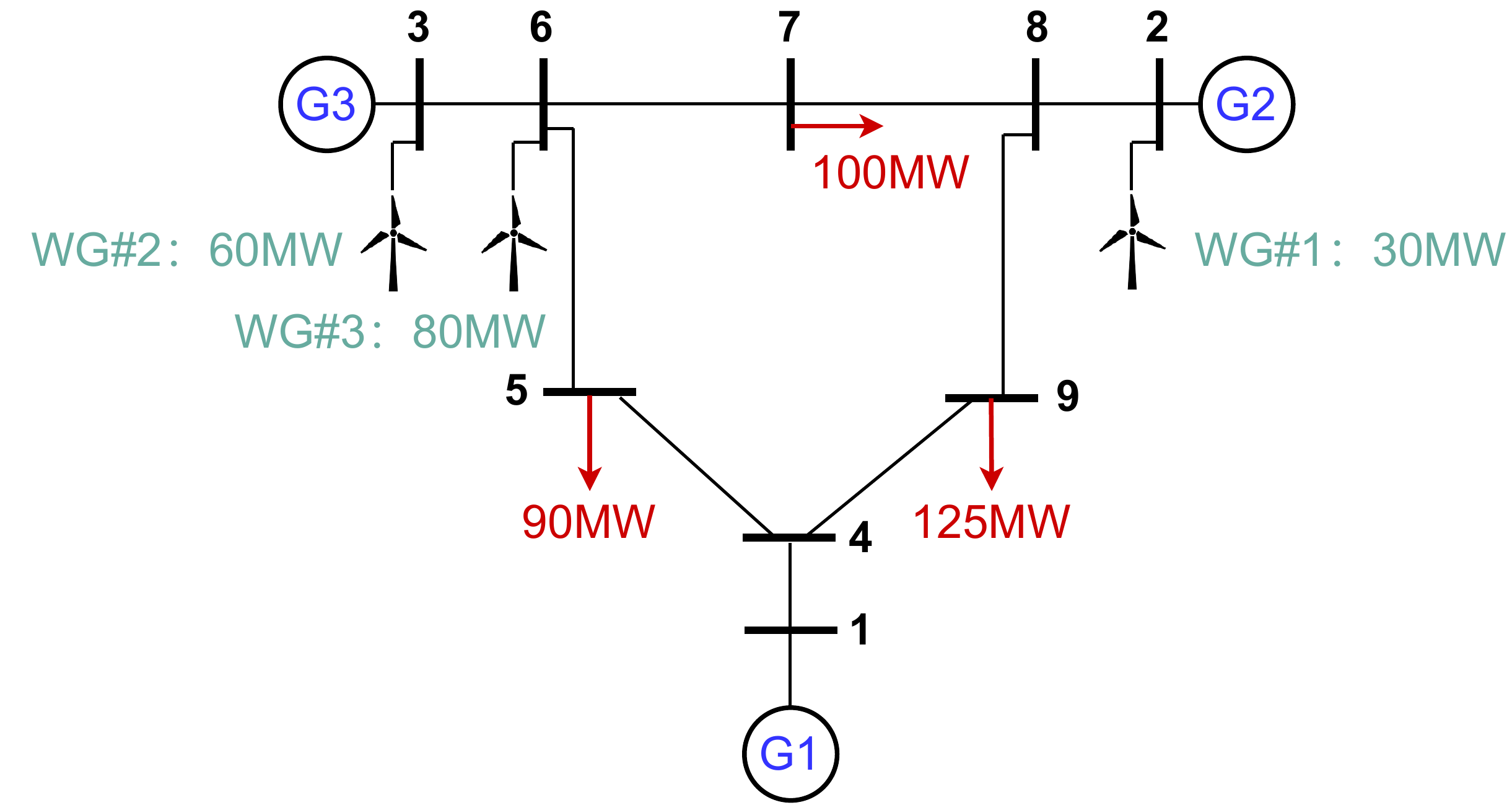}\vspace{-2mm}
	\caption{Structure of the test system. WG represents the wind generator (e.g. WG\#2:60MW denotes the wind generator with 60MW as its nominal output). Red arrows are power demands. } 
	\label{fig_system}
	\vspace{-2mm}
\end{figure}

The case study is based on the IEEE 9 bus system, as shown in Fig. \ref{fig_system}. The historical data of uncertainties are constructed based on the true wind farm SCADA data \cite{optis2019openoa}, and the corresponding probability density function of the dataset will be presented in Fig. \ref{fig_GMM}. The point number $N$ in  (\ref{eqn_inverse2}) for piece-wise linearization is set as 10. The parameters of generators can be found in Table \ref{tab_gene}. 

All numerical experiments are implemented on an Intel(R) 8700 3.20GHz CPU with 16 GB memory. We employ CVXPY to build our optimization problem and GUROBI to solve it.

\begin{table}[]
\small
\centering
\caption{Parameters of generators}
\begin{tabular}{ccccc} 
\hline
\multirow{2}{*}{Generator} & \multicolumn{2}{c}{{Generation bound (MW)}}          & \multirow{2}{*}{\makecell[c]{$c_1^i$ \\(\$/MWh$^2$)}} & \multirow{2}{*}{\makecell[c]{$c_2^i$ \\ (\$/MWh)}} \\ \cline{2-3}
                           & \multicolumn{1}{c}{Lower} & \multicolumn{1}{c}{Upper} &                                                    &                               \\ \hline
G1                         & 0                         & 250                       & 0.11                                               & 5                             \\
G2                         & 0                         & 300                         & 0.085                                               & 1.2                           \\
G3                         & 0                         & 270                       & 0.1225                                             & 1                             \\ \hline
\end{tabular}\label{tab_gene}
\vspace{-4mm}
\end{table}

To validate the superiority of the proposed model, we introduce two benchmarks, as follows:
\begin{enumerate}
	\item \textbf{B1}: CCP with Gaussian assumption.
	\item \textbf{B2}: Robust optimization model.
\end{enumerate}

\subsection{Effectiveness of GMM}
\begin{figure}
	\centering
	\vspace{-4mm}
	\subfigbottomskip=-6pt
	\subfigcapskip=-4pt
	\subfigure[]{\includegraphics[width=0.49\columnwidth]{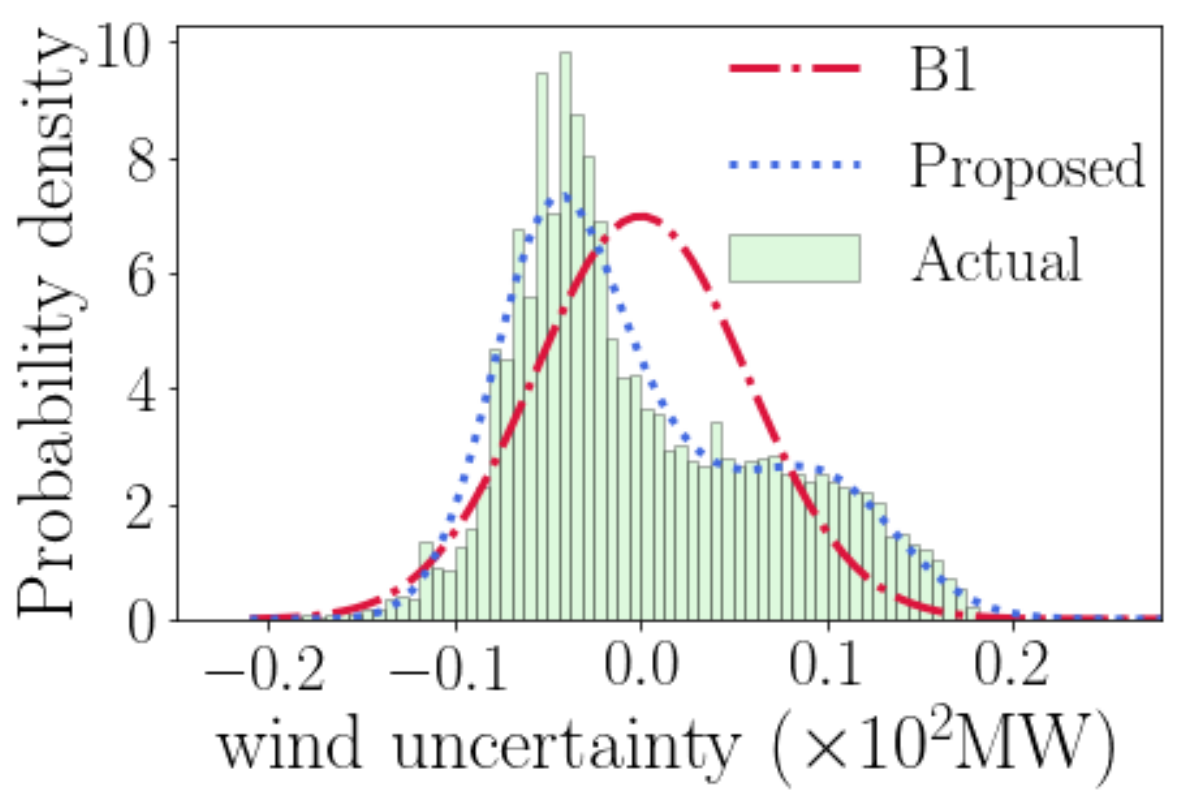}}
	\subfigure[]{\includegraphics[width=0.49\columnwidth]{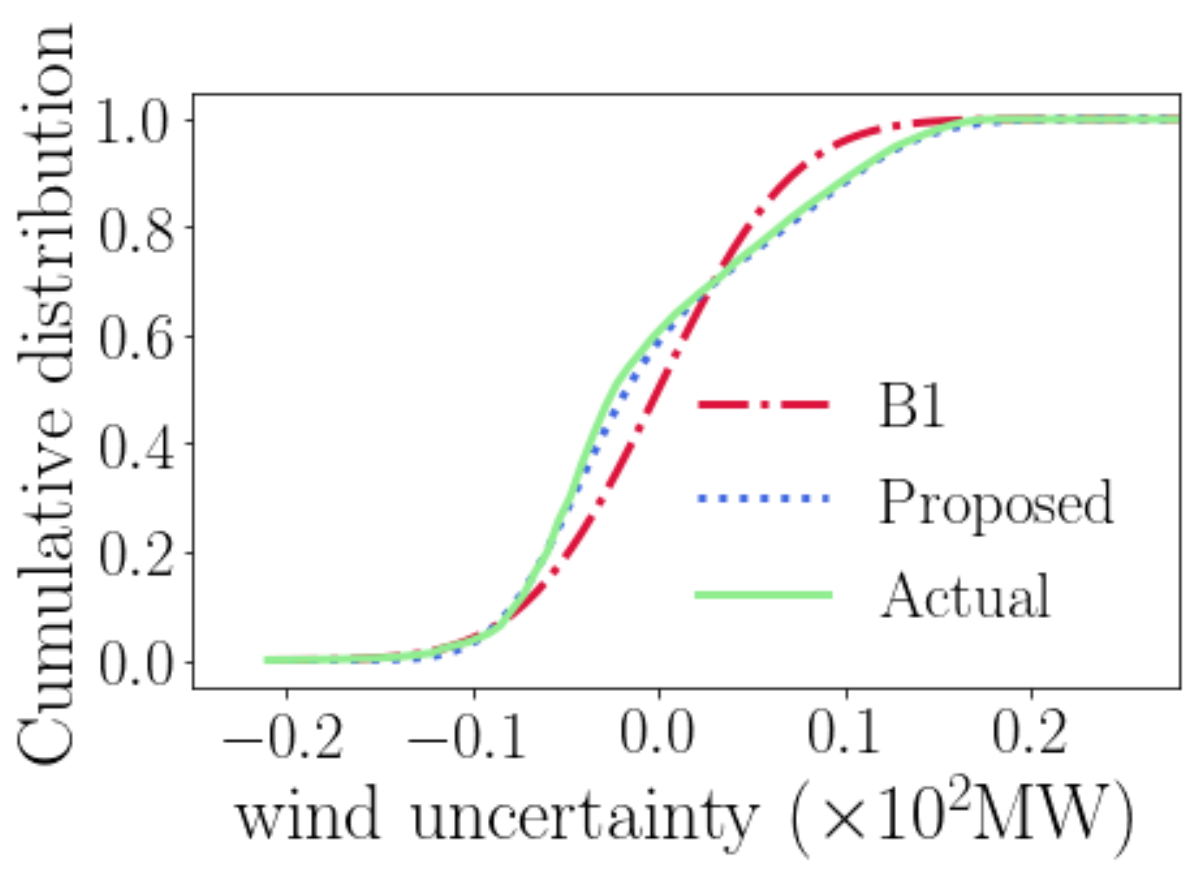}}
    \vspace{-4mm}
	\caption{Fitted PDFs and CDFs of the uncertainties from WG\#2. 
	}
	\label{fig_GMM}
	 		\vspace{-4mm}
\end{figure}

Fig. \ref{fig_GMM} demonstrates the fitted probability density functions (PDFs) and CDFs. The component number of GMM in the proposed model is 3. Obviously, the samples in the original dataset are non-Gaussian distributed. Thus, the actual PDF is far away from that of the Gaussian distribution. As a result, benchmark \textbf{B1} can not well fit the actual PDF since it is based on Gaussian assumption. On the contrary, the proposed GMM-based method can well capture the characteristics of the original irregular distribution, as shown in Fig. \ref{fig_GMM}(a). The distance from the fitted CDF to the actual one is also much smaller compared to \textbf{B1}, as shown in Fig. \ref{fig_GMM}(b). These results indicate the desirable fitting ability of GMM.

\subsection{Optimality and feasibility}
Table \ref{tab_results} summarizes the energy cost, solving time, and violation probability of different models. Among the three models, benchmark \textbf{B2} derives the most conservative results because it requires that the constraints should be satisfied with all realizations of uncertainties, which harms the energy efficiency of the system. The costs of \textbf{B1} are slightly lower than the proposed model. However, the results of violation probability are always higher than the given risk parameter $\epsilon$. That is to say, the reliability of \textbf{B1} can not ensure the probability requirement with a predetermined confidence level. Conversely, although the costs of the proposed model are higher than those in \textbf{B1}, the violation probabilities stay in the safe range, i.e, the probabilities are lower than $\epsilon$. In other words, the proposed GMM-based model can still ensure operation safety with a predetermined probability even with non-Gaussian distributed uncertainties on the left-hand side. These results show the high reliability of the proposed model.

Since the proposed model introduces binary variables for the linearization of the bilinear terms, the corresponding computational burden is larger than the rest of the models. However, as shown in Table \ref{tab_results}, in our case study, the solving times are around 1s, which is far less than the policy update interval in practice. Thus, the computational efficiency of the proposed model is also acceptable.

\begin{table}[]
\footnotesize
\centering
\caption{Results of all models under different risk parameters}
\begin{tabular}{cccccc}
\hline
\multicolumn{2}{c}{\textbf{Risk parameter $\epsilon$}}                            & \textbf{0.05} & \textbf{0.10} & \textbf{0.15}   & \textbf{0.20}  \\ \hline
\multirow{3}{*}{\makecell[c]{\textbf{Energy} \\\textbf{Cost (\$)}}}           & Proposed & 1050.54 & 1049.65 & 1049.14  &1048.57  \\ \cline{2-6} 
                                            & B1       & 1048.75 & 1047.99 & 1047.45  & 1047.04  \\ \cline{2-6}
                                            & B2       & \multicolumn{4}{c}{1065.38}    \\ \hline
\multirow{3}{*}{\makecell[c]{\textbf{Solving} \\\textbf{Time (s)}}}          & Proposed & 0.18 & 0.15 & 1.79 & 0.18 \\ \cline{2-6} 
                                            & B1       & 0.05 & 0.05 & 0.04  & 0.03  \\ \cline{2-6}
                                            & B2       & \multicolumn{4}{c}{0.03}    \\ \hline
\multirow{3}{*}{\makecell[c]{\textbf{Violation} \\\textbf{probability} \\ \textbf{(\%)}}} & Proposed & 4.01 & 7.69 & 10.51 &14.39  \\ \cline{2-6} 
                                            & B1       & 13.19 & 18.86 & 23.43  & 27.69  \\ \cline{2-6}
                                            & B2       & \multicolumn{4}{c}{0}    \\ \hline
\end{tabular} \label{tab_results}

\end{table}

\section{Conclusions} \label{sec_conclusion}
In this paper, a novel GMM-based CCP OPF model is proposed to handle the non-Gaussian distributed uncertainties from stochastic DG on the left-hand side. This model leverages GMM to fit the irregular CDF of the original samples. To ensure the tractability of the proposed model, a novel linearization technique is further developed to convert the original intractable terms into linear forms with binary variables that can be directly handled by off-the-shelf solvers. To verify the superiority of the proposed model, we conduct a case study based on the IEEE 9 bus system. Numerical results demonstrate that the proposed model can better capture the characteristics ofsamples' probability distribution. Simulation results also confirm that the proposed one can better ensure reliability than the CCP with Gaussian assumption. 


\footnotesize
\bibliographystyle{ieeetr}
\bibliography{ref}

\begin{thebibliography}{10}

\bibitem{IEA20202}
IEA, ``Renewables 2020.'' [OL].
\newblock \url{https://www.iea.org/reports/renewables-2020} Accessed November,
  2020.

\bibitem{8550809}
M.~{Islam}, M.~{Nadarajah}, and M.~J. {Hossain}, ``Short-term voltage stability
  enhancement in residential grid with high penetration of rooftop pv units,''
  {\em IEEE Trans. Sustain. Energy}, vol.~10, no.~4, pp.~2211--2222, 2019.

\bibitem{8016417}
A.~Lorca and X.~A. Sun, ``The adaptive robust multi-period alternating current
  optimal power flow problem,'' {\em IEEE Trans. Power Syst.}, vol.~33, no.~2,
  pp.~1993--2003, 2018.

\bibitem{bienstock2014chance}
D.~Bienstock, M.~Chertkov, and S.~Harnett, ``Chance-constrained optimal power
  flow: Risk-aware network control under uncertainty,'' {\em Siam Review},
  vol.~56, no.~3, pp.~461--495, 2014.

\bibitem{geng2019data}
X.~Geng and L.~Xie, ``Data-driven decision making in power systems with
  probabilistic guarantees: Theory and applications of chance-constrained
  optimization,'' {\em Annu Rev Control}, vol.~47, pp.~341--363, 2019.

\bibitem{8017474}
L.~Roald and G.~Andersson, ``Chance-constrained ac optimal power flow:
  Reformulations and efficient algorithms,'' {\em IEEE Trans. Power Syst.},
  vol.~33, no.~3, pp.~2906--2918, 2018.

\bibitem{pena2020dc}
A.~Pena-Ordieres, D.~K. Molzahn, L.~A. Roald, and A.~W{\"a}chter, ``Dc optimal
  power flow with joint chance constraints,'' {\em IEEE Trans. Power Syst.},
  vol.~36, no.~1, pp.~147--158, 2020.

\bibitem{9535415}
G.~Chen, B.~Yan, H.~Zhang, D.~Zhang, and Y.~Song, ``Time-efficient strategic
  power dispatch for district cooling systems considering the spatial-temporal
  evolution of cooling load uncertainties,'' {\em CSEE Journal of Power and
  Energy Systems}, pp.~1--11, 2021.

\bibitem{7307226}
D.~Ke, C.~Y. Chung, and Y.~Sun, ``A novel probabilistic optimal power flow
  model with uncertain wind power generation described by customized gaussian
  mixture model,'' {\em IEEE Trans. Sustain. Energy}, vol.~7, no.~1,
  pp.~200--212, 2016.

\bibitem{8936474}
Y.~Yang, W.~Wu, B.~Wang, and M.~Li, ``Analytical reformulation for stochastic
  unit commitment considering wind power uncertainty with gaussian mixture
  model,'' {\em IEEE Trans. Power Syst.}, vol.~35, no.~4, pp.~2769--2782, 2020.

\bibitem{8772186}
W.~Sun, M.~Zamani, M.~R. Hesamzadeh, and H.-T. Zhang, ``Data-driven
  probabilistic optimal power flow with nonparametric bayesian modeling and
  inference,'' {\em IEEE Trans. Smart Grid}, vol.~11, no.~2, pp.~1077--1090,
  2020.

\bibitem{9376652}
J.~Wang, C.~Wang, Y.~Liang, T.~Bi, M.~Shafie-khah, and J.~P.~S. Catalao,
  ``Data-driven chance-constrained optimal gas-power flow calculation: A
  bayesian nonparametric approach,'' {\em IEEE Trans. Power Syst.}, pp.~1--1,
  2021.

\bibitem{7973099}
W.~Xie and S.~Ahmed, ``Distributionally robust chance constrained optimal power
  flow with renewables: A conic reformulation,'' {\em IEEE Trans. Power Syst.},
  vol.~33, no.~2, pp.~1860--1867, 2018.

\bibitem{goodfellow2016deep}
I.~Goodfellow, Y.~Bengio, A.~Courville, and Y.~Bengio, {\em Deep learning},
  vol.~1.
\newblock MIT press Cambridge, 2016.

\bibitem{hu2018chance}
Z.~Hu, W.~Sun, and S.~Zhu, ``Chance constrained programs with mixture
  distributions,'' 2018.

\bibitem{optis2019openoa}
M.~Optis, J.~Perr-Sauer, C.~Philips, A.~E. Craig, J.~C. Lee, T.~Kemper,
  S.~Sheng, E.~Simley, L.~Williams, M.~Lunacek, {\em et~al.}, ``Openoa: An
  open-source code base for operational analysis of wind power plants,'' {\em
  Wind Energy Science Discussions}, pp.~1--14, 2019.

\end{thebibliography}
\end{document}